\documentclass[10pt,a4paper,fleqn]{amsart}
\usepackage{amssymb,amsmath,latexsym}
\usepackage[varg]{pxfonts}
\usepackage[english]{babel}
\def\constr#1^#2{\mathrel{\mathop{\kern 0pt#1}\limits^{#2}}}
\def\build#1_#2{\mathrel{\mathop{\kern 0pt#1}\limits_{#2}}}

\theoremstyle{plain}

\newtheorem{Theorem}{Theorem}[section]

\theoremstyle{Corollary}
\newtheorem{Corollary}{Corollary}[section]

\newtheorem{Lemma}{Lemma}[section]
\theoremstyle{Proposition}

\theoremstyle{definition}
\newtheorem{Definition}{Definition}[section]

\newtheorem{remark}{Remark}[section]

\theoremstyle{Abstract}

\numberwithin{equation}{section}

\newcommand{\pv}{\par \vskip0.2cm}

\newcommand{\hf}{\hfill $\diamondsuit$ }    

\newenvironment{pf}{\noindent{\bf Proof.}\enspace}{
\hfill }

\begin{document}
\title[]{Fixed Point Results for Multivalued Mappings with Closed Graphs in Generalized Banach Spaces}
\author[]{ Khaled Ben Amara$^{(1)},$  Aref Jeribi $^{(1)},$  Najib Kaddachi $^{(2)}$ and  Zahra Laouar$^{(1)}$}
\date{}
\maketitle
\begin{center}
 Department of Mathematics. Faculty of Sciences of Sfax.
University of Sfax.\\ Road Soukra Km $3.5 B.P. 1171, 3000,$ Sfax,
Tunisia.\pv
\end{center}

\vskip0.3cm

\centerline {e-mail: $^{(1)}$ Aref.Jeribi$$@$$fss.rnu.tn,   \
$^{(2)}$ najibkadachi$$@$$gmail.com   }
\vskip0.7cm
\vskip0.2cm

\maketitle

\par\vskip0.1cm

\textbf{Abstract.} In this paper, by establishing a new characterization of the notion of upper semi-continuity of multivalued mappings in generalized Banach spaces, we prove some Perov type fixed point theorems for multivalued mappings with closed graphs. Moreover, we derive some Krasnoselskii's fixe point results for  multivalued mappings in generalized Banach spaces. Our results are applied to a large class of nonlinear inclusions systems.

\vskip0.7cm

\section{\textbf{Introduction}}
 Fixed point theory for
multi-valued mappings plays an important role in the study of many problems in several fields of applied sciences, see \textup{\cite{BA2015,Dhage,JK2015,Kaddachi2016}} and the references therein contained. Different directions were followed  in order to develop this theory.
 \begin{Theorem}\cite{algebra}\label{usc}
Let $S$ be a nonempty, compact and convex subset of a generalized Banach space  $X$ and let $F : S \longrightarrow   \mathcal{P}_{cp,cv}(S)$ be an upper semi-continuous  multi-valued operator. Then  $F$ has at least one fixed point, that is, there exists $x \in S$ such that
$x \in  F(x).$
\end{Theorem}
Later, in \cite{measureg}, J. R. Graef, J. Henderson and A. Ouahab have introduced the concept of measure of noncompactness
in vector-valued Banach spaces, in order to obtain a Sadovkii's  fixed point theorem for singlevalued mappings.

 In this paper, we establish a characterization of the upper semi-continuity of a multivalued mapping in generalized Banach space by the closeness of its graph.
Moreover, we prove some new  fixed point results of Perov's type for closed multi-valued operator. This allows us to derive some Krasnoselskii's fixed point theorems,  without upper semi-continuity condition.

 we give in the next section some 
preliminaries which will be used throughout this paper.  In section $3,$ by invoking the Hausdorff measure of noncompactness, we extend some classical results due to  Nussbaum \cite{propo}, and represent the relationship between  $k$-Lipschitzian and $k$-set-Lipschitzian mappings defined on a classical metric spaces to the generalized metric spaces. 
\section{\textbf{Preliminaries}}
\noindent At the beginning, we shall introduce some definitions and give preliminary results which will be needed in the sequel.
 If, $x:=(x_1, \ldots ,x_n), y:= (y_1, \ldots , y_n) \in \mathbb{R}^n,$
 by $x \leq y$ we mean $x_i \leq y_i$ for all $i = 1, \ldots, n.$ Given $\alpha \in \mathbb{R},$ by $x \leq \alpha$ we mean $x_i \leq \alpha$ for
all $i = 1, \ldots, n.$ For all $x:=(x_1, \ldots ,x_n), y:= (y_1, \ldots , y_n) \in \mathbb{R}^n,$ we put $$|x| = (|x_1|, \ldots , |x_n|), \hbox{ and }\max(x, y) = (\max(x_1, y_1), \ldots, \max(x_n, y_n)).$$
\begin{Definition}\cite{perov1}
Let $X$ be a nonempty set. A mapping $d:X\times X\longrightarrow \mathbb{R}^{n}$ is called a vector-valued metric on $X,$ if the following statements are satisfied\\
\noindent{$(i)$} $d(x,y)\geq 0$ for all $x,y \in X,$ if $d(x,y)=0$ then $x=y,$\\
\noindent{$(ii)$} $d(x,y)=d(y,x)$ for all $x,y \in X,$\\
\noindent{$(iii)$} $d(x,y)\leq d(x,z)+d(z,y)$ for all $x,y,z \in X.$\hf
\end{Definition}
\noindent We call the pair $(X, d)$ a generalized metric space with
$$d(x,y):=\left(%
\begin{array}{ccc}
  d_1(x,y)  \\
  \vdots \\
   d_n(x,y) \\
\end{array}%
\right).$$
Notice that $d$ is a generalized metric on $X,$ if and only, if $d_i, i=1,\ldots,n,$ are metrics on $X.$
\begin{Definition}\cite{algebra,PetreIoanRadu,IOAN}
Let $X$ be a vector space on $\mathbb{K}$\ \  $(\mathbb{K}=\mathbb{R} \text{   or   } \mathbb{K}=\mathbb{C}).$ By a vector-valued norm on $X.$ We mean a map $\|\cdot\|:X\longrightarrow \mathbb{R}^{n}_{+}$ with the following properties:\\
\noindent${(i)}$ $\|x\|\geq 0$ for all  $x \in X,$ if $\|x\|=0$  then  $x=(0,...,0),$\\
\noindent${(ii)}$ $\|\lambda x\|= |\lambda|\|x\|$ for all $x \in X$ and $\lambda \in \mathbb{K},$\\
\noindent${(iii)}$ $\|x+y\|\leq \|x\|+\|y\|,$ for all $x,y \in X.$\hf
\end{Definition}
\noindent The pair $(X,\|\cdot\|)$ is called a generalized normed space. If the generalized metric
generated by  $ \|\cdot\|$  is complete then the space $(X, \|\cdot\|)$ is called a generalized Banach space. Moreover,  if $\|\cdot\|_i, i=1, \ldots,n,$ are norms on $E,$ then $(X,\|\cdot\|),$ where $\|\cdot\|:=\left(%
\begin{array}{ccc}
  \|\cdot\|_1  \\
  \vdots \\
  \|\cdot\|_n \\
\end{array}%
\right),$ is a generalized Banach space, see \cite[P. 135]{measureg}.\\
Now, for $r=(r_1,\ldots,r_n) \in \mathbb{R}^{n}_{+},$ we will denote by
$$B(x_0,r)=\{x \in X,\ \ \|x-x_0\|< r\}\ \ \text{  the open ball centered in } x_0 \ \ \text{with radius}\ \ r,$$ and by
$$\overline{B(x_0,r)}=\{x \in X, \ \ \|x-x_0\|\leq r\}\ \ \text{  the closed ball centred in } x_0 \ \ \text{with radius}\ \ r.$$

\noindent Let $(X,d)$ be a generalized metric space, we denote by

\noindent Let $S$ be a nonempty subset of $X$ and $F : S \to \mathcal{P}(X)$ be a multi-valued mapping. By $Gr(F),$ we denote the graph of $F$ which is defined by the formula
\begin{equation*}
Gr(F)=\left\{(x,y) \in  S\times X:\, y \in T(x)\right\}.
\end{equation*}
Moreover, for all subset $M$ of $X,$ we
put \begin{equation*}F^{-1}(M) = \left\{x \in S :\, F(x) \cap M\neq \emptyset \right\} \ \ \hbox{and} \ \ F^{+}(M) = \left\{x \in S :\, F(x) \subset M \right\}.\end{equation*}

\begin{Definition}\cite{IOAN}
 Let $F : S \longrightarrow \mathcal{P} (X ).$  We say that
 \\
 \noindent $(i)$ $F$ is upper semi-continuous (u.s.c) if $F^{+}(V)$ is open for any open $V $ of $X,$\\
 \noindent $(ii)$ $F$ has a closed graph if for every sequence $\{x_n, n \in \mathbb{N}\}\subset S$ with $x_n \rightarrow x \in S$ and for every sequence $\{y_n, n \in \mathbb{N}\}$ with $y_n \in F(x_n),$ such that $y_n \rightarrow y,$ then $y \in F(x).$\hf
 \end{Definition}
\noindent Let $(X, d_*)$ be a metric space, we denote by $H_{d_*}$ the Hausdorff pseudo-metric
distance on $\mathcal{P}(X)$ defined as
$$H_{d_*}: \mathcal{P}(X)\times \mathcal{P}(X)\longrightarrow \mathbb{R}_{+} \cup \{+\infty\}, \ \
H_{d_*}(A,B)=\max  \left\{\sup_{a \in A} D_*(a,B);\sup_{b \in B} D_*(A,b)  \right\} $$
where $D_*(A,b)=\inf_{a \in A}d_*(a,b) $ and  $ D_*(a,B)=\inf_{b \in B}d_*(a,b).$\\
Then $(\mathcal{P}_{bd,cl}(X), H_{d^*})$
is a metric space and $(\mathcal{P}_{cl}(X), H_{d^*})$ is a generalized metric space. In particular, $H_{d^*}$
satisfies the triangle inequality, see\cite{OUAHAB,IOAN}.
Now, we consider the generalized Hausdorff pseudo-metric distance
$$H_{d}: \mathcal{P}(X)\times \mathcal{P}(X)\longrightarrow \mathbb{R}^{n}_{+} \cup \{+\infty\},$$
defined by
$$H_{d}(A,B):=\left(%
\begin{array}{ccc}
  H_{d_1}(A,B)  \\
  \vdots \\
  H_{d_n}(A,B) \\
\end{array}%
\right),$$
where $H_{d_i},i=1,\ldots,n$ are Hausdorff pseudo-metric distance.

\begin{Definition}\cite{algebra,IOAN}
Let $(X, d)$ be a generalized metric space. A multivalued mapping $ F:X\longrightarrow \mathcal{P}_{cl}(X)$ is said to be:\\
$(i)$ Lipschizian if there exists a matrix $M \in  \mathcal{M}_{n \times n}(\mathbb{R}_+)$ such that $$H_d(Fx,Fy)\leq Md(x,y),\ \ \text{ for all }  x, y \\  \text { in } X. $$
$(ii)$ Contraction in the sense of Nadler if there exists a matrix $M \in  \mathcal{M}_{n \times n}(\mathbb{R}_+)$ convergent to zero such that
$$H_d(Fx,Fy)\leq Md(x,y),\ \ \text{ for all }  x, y \\  \text { in } X. $$\hf
\end{Definition}
\begin{Definition}\cite{measureg}
Let $X$ be a generalized Banach space and $(\mathcal{A},\leq)$ be a partially ordered
set.  A map $\mu : \mathcal{P}(X) \longrightarrow \mathcal{A} \times \mathcal{A}\times \ldots \times \mathcal{A} $ is called a generalized measure of noncompactness on $ X$, if
$$ \mu(\overline{co }(V))=\mu(V), \ \ \text { for every  } V \in \mathcal{P}(X)$$
where $$\mu(V):=\left(%
\begin{array}{ccc}
  \mu_1(V)  \\
  \vdots \\
  \mu_n(V)  \\
\end{array}%
\right).$$\hf
\end{Definition}
\begin{Definition} \cite{measureg}\label{Defmg}
A measure of noncompactness $\mu$ is called:\\
\noindent {$(i)$} Monotone if $V_1\subseteq V_2$ then   $\mu(V_1)\leq \mu(V_2)$ for all $V_1,V_2 \in \mathcal{P}(X).$\\
\noindent {$(ii)$} Nonsingular if $\mu(\{a\}\cup V )=\mu(V)  $ for every $a \in X$ and $V \in \mathcal{P}(X).$\\
\noindent {$(iii)$} Real if $ \mathcal{A} = \overline{\mathbb{R}}_+$ and $\mu(V)< \infty$ for every $i = 1,\ldots, n$ and every bounded $V.$\\
\noindent {$(iv)$} Semi-additive if $\mu(V_1 \cup V_2) = \max(\mu(V_1), \mu(V_2))$ for every $V_1,V_2 \in \mathcal{P}(X).  $ \\
\noindent {$(v)$} Lower-additive if $\mu(V_1+V_2)\leq \mu (V_1)+\mu(V_2),$ for all $V_1,V_2 \in \mathcal{P}(X).$\\
\noindent {$(vi)$} Regular if the condition $\mu(V_1) = 0_{\mathbb{R}^n}$ is equivalent to the relative compactness of $V.$\hf
  \end{Definition}
  \begin{Definition}\cite{measureg}
 Let $X$ be a generalized metric space. A multivalued map $ F:X\longrightarrow \mathcal{P}_{cl}(X)$ is said to be:\\
$(i)$ Lipschizian if there exists a matrix $M \in  \mathcal{M}_{n \times n}(\mathbb{R}_+)$ such that $$ \mu(F(V)) \leq  M\mu(V), \ \  \text{  for all } V \in  \mathcal{P}(X).$$
$(ii)$ Contraction in the sense of Nadler if there exists a matrix $M \in  \mathcal{M}_{n \times n}(\mathbb{R}_+)$ convergent to zero such that
$$ \mu(F(V)) \leq  M\mu(V), \ \  \text{  for all } V \in  \mathcal{P}(X).$$ \hf
 \end{Definition}
\begin{remark}\label{remarkgéfiny}
From the  Definition $\ref{Defmg},$ if $\mu$ is a real genaralised measure of non compactness and assume that  $F: X \longrightarrow \mathcal{P} (Y )$ is $M$-contraction multivalued map with respect to $\mu,$ then $\mu(F(V))<\infty$ for every bounded $V $ of $\mathcal{P}(X)$ \hf
\end{remark}
\section{\textbf{Perov's fixed point results for multivalued mappings with closed graphs}}
\noindent In this section, we are going to discuss some fixed point theorems  involving nonlinear multivalued mappings in generalized Banach space by using the  concept of the measure of non compactness. Notice that Deimling in  \cite{Deimling} gave a characterization  of upper semi-continuous for a multivalued mappings in Banach space  in terms of the closure of its graph.
 Before stating the main results, we  state  Lemmas useful in the sequel which are an extension  of the Deimling's results to the generalized metric space.
\begin{Lemma} \label{pro} Let $S$ be a nonempty subset of a   generalized  Banach space $E$ and let $Q: S \subset E \longrightarrow \mathcal{P}_{cl}(E)$ be a  multi-valued operator  such that $Q^{-1}(V)$ is closed, for any closed subset  $V$ of $E.$ Then $Q$ is upper semi-continuous.
\end{Lemma}
\begin{pf}
Let $E_i = E; i = 1, \ldots, n.$ Consider
$\displaystyle\prod_{i=1}^{n}E_i$ with $$ \|(x_1,\ldots,x_n)\|_{\ast}=  \sum_{i=1}^{n}\|x\|_i \ \   \text{  for each } x \in E.$$
The diagonal space of $\displaystyle\prod_{i=1}^{n}E_i$ defined by
$$ \widetilde{ E}= \big\{ (x,\ldots,x)  \in \displaystyle \prod_{i=1}^{n}E_i : x \in E, i=1,\ldots,n \big \}. $$
Let us define the mapping $h : E\longrightarrow \widetilde{ E}$ by the formula
$
h(x)=(x,\ldots,x),
$
where $$ \widetilde{ E}= \big\{ (x,\ldots,x)  \in \displaystyle \prod_{i=1}^{n}E_i \ \  \text{ such that } x \in E_i=E; i=1,\ldots,n \big \}. $$
Notice that $\widetilde{E}$ is a Banach space endowed with the norm $\|\cdot\|_{\ast}$ defined by
$$ \|(x,\ldots,x)\|_{\ast}=  \sum_{i=1}^{n}\|x\|_i \ \   \text{  for each } x \in E,$$
where $\|\cdot\|_i,i=1,\ldots,n$ are norms. Let $R :=h\circ Q \circ h^{-1}:h(S)\subset \widetilde{ E} \longrightarrow \mathcal{P}_{cl}(\widetilde{ E}). $
The use of our assumption  as well as  Lemma $7.12$ in   \cite{measureg} allows us to say that  $Q^{-1}(h^{-1}(A))$ is closed, for any closed subset $A$ of $\widetilde{E}.$
 Using \cite[page.3]{cross}, we get $R^{-1}(A)$
is closed in $h(S).$ As a
result, $ R$  is upper semi-continuous on $h(S)$ in light of Proposition $24.1$ in \cite{Deimling}. Now, we prove that $Q$ is upper semi-continuous on $E.$ To see this, let $U$ be an open on $E,$ then $h(U)$ is an open on $\widetilde{E}.$
Let us  consider
\begin{align*}
 R^+(h(U))& = \big\{ y \in h(S),  R(y) \subset h(U)\big\}
            = \big\{h(x), \ \ x \in S  \ \  R(h(x))\subset  h(U)  \big\}\\
            &=  h \left(\big\{x \in S, \ \ R(h(x))\subset  h(U) \big\}\right).
 \end{align*}
Since  $ R$  is upper semi-continuous on $h(S)$, it follows that

\begin{align}\label{h(u)open}
R^+(h(U))=  h \left(\left\{x \in S, \ \ h(Q(x))\subset  h(U) \right\}\right)
               =  h \left(\left\{x \in S, \ \ Q(x)\subset  U \right\}\right)
               = h \left( Q^+(U)\right)
\end{align}
 is an open subset of $\widetilde{ E}$  and consequently $Q^+(U)$ is open in $E$ in view of proposition $24.1$ in \cite{Deimling}.
\end{pf}

\bigskip

\begin{Lemma}\label{uscclosedgraph1}Let $S$ be a nonempty and closed subset of a generalized Banach space $E$ and
let  $Q: S \subset E\longrightarrow \mathcal{P}_{cl}(E)$ be a multi-valued operator.
 If $Q$ is upper semi-continuous, then $Q$ has a closed graph.\hf
\end{Lemma}

\begin{pf}
Let $U$ be an open in $\widetilde{E}.$ Then $h^{-1}(U)$ as well as $Q^{+}(h^{-1}(U))$ is open in $E.$
 Using the equality $(\ref{h(u)open})$ combined with the fact that $h$  is
a homeomorphism, we obtain that $R^{+}(U)$ is an open on $\widetilde{E},$ where $R$ is defined in Lemma \ref{pro}, or equivalently $R$ is upper semi-continuous. From proposition $24.1$ in \cite{Deimling}, it follows that $R$ has a closed graph. Next, we claim that $Q$ has a closed graph.  To
see this, consider $\{x_n, \ \ n \in \mathbb{N} \} $ as a convergent sequence in $S$ to a point $x$ and let $\{y_n, \ \ n \in \mathbb{N} \}$ be a  converging sequence in $Q(x_n)$ to a point $y.$ Making use of   Lemma $7.12$ in   \cite{measureg}, together with the assumption on $Q,$ enables us to have
 $h(x_n) \rightarrow h(x)$ and $\{h(y_n), \ \ n \in \mathbb{N} \} $ is a  converging sequence in  $ R(h(x_n)))$  to $h(y).$
From the above discussion, it is easy
to see that  $ h(y) \in  R(h(x)))$ or equivalently   $y \in Qx$ and the claim is proved.
\end{pf}

\bigskip

\begin{Lemma}\label{uscclosedgraph}Let $S$ be a nonempty and closed subset of a generalized Banach space $E$ and
let  $Q: S  \longrightarrow \mathcal{P}_{cl}(E)$ be a multi-valued operator such that $Q(S)$
 is relatively compact. Then $Q$ has a closed graph  if and only if  $Q$ is upper semi-continuous. \hf
\end{Lemma}

\begin{pf}
In view of Lemma \ref{uscclosedgraph1}, we need only to prove necessity. Firstly, we claim that   $R$ has a closed graph, where $R$ is defined in Lemma \ref{pro}. To see this, let $\{\xi_n,\ \ n \in \mathbb{N}\}$ be a sequence in $ h(S),$ with  $\xi_n\rightarrow \xi \in h(S)$ and let
$\nu_n \in R(\xi_n), n \in \mathbb{N},$ with $\nu_n \rightarrow \nu.$
 If we take the sequences $\{x_n,\ \ n \in \mathbb{N}\}$ and $\{y_n,\ \ n \in \mathbb{N}\}$ defined
by
 $$x_n=h^{-1}(\xi_n) \ \ \text { and  } \ \  y_n=h^{-1}(\nu_n), \ \ n \in \mathbb{N},$$
then $x_n \rightarrow h^{-1}(\xi)=x \in E$ and
$y_n \in Q(x_n) $ with $y_n \rightarrow h^{-1}(\nu )=y. $ Since $Q$ has a closed graph $E,$ then
$y \in Q(x).$ This implies that $\nu  \in R(\xi),$ and the claim is proved.
  Our next task is to prove that $Q^{-1}(V)$ is closed for all closed $V$ in $E.$ For this purpose, let  $V$ be a closed subset of $E$ and let $(x_n)_n \subset Q^{-1}(V)$ such that $x_n \rightarrow  x.$
Then, there exists a sequence $\{y_n, \ \ n \in \mathbb{N}\}$ such that
\begin{eqnarray}\label{intrsection}
y_n \in Q(x_n)\cap V, \ \  \text{ for all  }  n \in \mathbb{N}.
\end{eqnarray}
Since $Q(S)$ is relatively compact, we can extract a subsequence  $(y_{n_k})_k$ of $\{y_{n}, \ \ n \in  \mathbb{N} \}$ such that
$$y_{n_k}\rightarrow y$$
 From (\ref{intrsection}), we get
$$ h(y_{n_k})  \in R(h(x_{n_k})), \ \ k \in \mathbb{N}.$$ Taking into account the fact that $R$ has a closed graph, we infer that $$h(y) \in R(h(x)).$$ This implies that
$$y \in  Q (x)  \cap V,$$ and consequently $$x \in Q^{-1}(V).$$ This achieves the claim and completes the proof, using Lemma \ref{pro}.
\end{pf}

\bigskip

\noindent  In the remained of this section, we  consider a monotone, regular, real generalized measure of non compactness $\mu: \mathcal{P}(X)\longrightarrow  \overline{\mathbb{R}}_+^{n}$ defined by
  \begin{eqnarray}\label{murealg}
   \mu(V):=\left(%
\begin{array}{ccc}
  \mu_1(V)  \\
  \vdots \\
  \mu_n(V)  \\
\end{array}%
\right).
\end{eqnarray}

\noindent In order to weakening the stability condition $(iii)$ of the above Theorem,  we establish the following result.
\begin{Theorem}\label{krasinve}
Let  $S$ be a nonempty, bounded, closed and convex subset of a generalized Banach space $E.$  Let
$F: E \longrightarrow   \mathcal{P}_{cl,cv}(E)$ and $G: S \longrightarrow   \mathcal{P}(E)$ be two multivalued mappings with closed graphs such that:\\
\noindent{$(i)$} There exists a continuous, nondecreasing function $\varphi: \mathbb{R}^n_+ \longrightarrow \mathbb{R}^n_+$ with $\displaystyle\lim_{m\rightarrow \infty}\varphi^m(t)= 0_{\mathbb{R}^n}$ for each $t \in   \mathbb{R}^n_+,$ such that
$$\mu(F(V))\leq \varphi(\mu(V))\hbox{ for all } V\in \mathcal{P}_{bd}(X),$$
\noindent{$(ii)$}  $F(E)$ is bounded and $G(S)$ is relatively compact,\\
\noindent{$(iii)$} $(I-F)^{-1}Gx$ is convex, for all $x \in S,$\\
\noindent {$(iv)$} $x \in F(x)+G(y), \ \ y \in S \Longrightarrow  x \in S.$\\
Then $G+F$ has, at least, one fixed point in $S.$\hf
\end{Theorem}
\begin{pf}
 Taking into account the assumptions $ (iii)$ and $(iv),$ we can define a  multivalued mapping $Q : S \longrightarrow\mathcal{P}_{cv}(S)$ by $$Qx=(I-F)^{-1}G(x).$$   Let us claim that $Q$ has closed values on $S.$ For this purpose, let $x \in S$ and let $\{y_n, \ \ n \in \mathbb{N} \}\subset Qx$ such that $y_n  \rightarrow y.$ Since $y_n\in Q(x),$ then $$(I-F)(y_n)\cap G(x) \neq \emptyset, n\in \mathbb{N}.$$ This implies that there exists $u_n \in F(y_n)$ and $v_n \in G(x)$ such that
 \begin{eqnarray}\label{a}y_n =u_n +v_n, n\in \mathbb{N}.\end{eqnarray}
 The fact that $G(S)$ is relatively compact allows $\{v_n, \ \  n \in \mathbb{N}\}$ to have  a renamed subsequence such that $$v_n\rightarrow v.$$ On the other hand, since $G$ has a closed graph, we get $v\in G(x).$
  Taking into account both  equality (\ref{a}) and the fact that $F$ has a closed graph, we infer that $$y-v \in F(y).$$ Hence, we have $$(I-F)(y) \cap G(x)\neq \emptyset,$$  which proves the claim.
Our next task is to prove that $Q(S)$ is relatively compact.
 In view of inclusion
 $$Q(S)\subseteq F(Q(S))+G(S)$$
 we have
 \begin{eqnarray*}
 \mu \left(Q(S)\right) \leq \varphi(\mu\left((Q(S)\right)).
 \end{eqnarray*}
Since $\varphi $ is nondecreasing, we get
$$ \mu(Q(S))\leq \varphi(\mu(Q(S))) \leq  \varphi^2(Q(S))) \leq \ldots \leq \varphi^m(\mu(Q(S))).$$
Letting $m \rightarrow \infty,$ and using the regularity of the generalized measure of noncompactness $\mu,$ we deduce the relatively compactness of $Q(S).$
 Since $Q(S)\subset S$ and $Q$ has closed values on $S,$ then we can consider the multi-valued mapping
 $Q$ by $$\left\{
                           \begin{array}{ll}
                              Q: \overline{co}(Q(S))\longrightarrow \mathcal{P}_{cp, cv}(\overline{co}(Q(S)))\\\\
                             x\mapsto \displaystyle (I-F)^{-1}G(x)
                           \end{array}
                         \right.
$$
In order to achieve the proof, we will apply  Theorem \ref{usc}.  Hence,
  in view of lemma \ref{uscclosedgraph}, it is sufficient to prove that $Q$ has a closed graph on $\overline{co}(Q(S)).$ To see this, let  $\{x_n, \ \ n \in \mathbb{N}\}$  be a sequence in $\overline{co}(Q(S))$ with  $x_n \rightarrow x$ and let  $y_n \in Q(x_n)$  with  $y_n \rightarrow y.$ This enables us to have $$y_n \in F(y_n)+G(x_n) , n \in \mathbb{N},$$ and consequently  there exist  $u_n \in F(y_n)$ and $v_n \in G(x_n),$ such that
 $$y_n=u_n+v_n, n\in \mathbb{N}.$$
 Therefore, from the assumption  $(ii),$ and in view of inclusion $\{v_n, \ \ n \in \mathbb{N}\} \subset G(\{x_n, \ \ n \in \mathbb{N}\}),$ we deduce that
 $\{v_n, \ \ n \in \mathbb{N}\}$ has a renamed subsequence which converges weakly to a point $v$ in $E.$
  Since $ u_{n_k} \rightarrow y-v,$ so by using the fact that
   $F$ and $G$ have closed graph, we get
 $$v \in (I-F)y \cap G(x).$$  This prove the claim and completes the proof by using Theorem \ref{usc}.
\end{pf}
\bigskip

\section{ \textbf{Multivalued fixed point results in terms of the Hausdorff measure of non compactness}}
This section is devoted to develop some fixed point results of Perov type for multivalued mappings, in terms of the Hausdorff measure of noncompactness.
\noindent Thus, we start by introducing a real generalized measure of non compactness $\mu: \mathcal{P}(E)\longrightarrow  \overline{\mathbb{R}}_+^{n},$ which will be used in this section,  by
  \begin{eqnarray}\label{murealg}
   \mu(V):=\left(%
\begin{array}{ccc}
  \mu_1(V)  \\
  \vdots \\
  \mu_n(V)  \\
\end{array}%
\right).
\end{eqnarray}
where the maps
  $$\mu_i: \mathcal{P}(E)\longrightarrow \overline{\mathbb{R}}_+, i=1,...,n$$ are defined by
\begin{eqnarray}\label{mureal}
\mu_i(V):=\left\{
 \begin{array}{ccc}
\chi_i(V) & \text {    if  }    & V \text{     bounded    }\\
  +\infty  &  \text {    if   }    & V \text{     unbounded }
\end{array}
 \right.
\end{eqnarray}
  where $\chi_i,\ i=1,\ldots,n,$ is the Hausdorff measure of noncompactness in the Banach spaces $(E,\|\cdot\|_i),$ see \cite{measure}.
 It is clear that $\mu$  satisfies the properties  $(i)-(vi)$.

\bigskip

\noindent In  \cite{propo}, Roger D. Nussbaum has proved that every $k$-contraction mappings defined on metric spaces is  $k$-set-contraction. This result was extended to the case of multivalued mappings, see \cite{Deimling}.
 The following Lemma present an extension of these results to the case of  multivalued mappings defined on  generalized Banach space.
\begin{Lemma}\label{sofgeneralisé}
 Let $E$ be a generalized Banach space. If the multivalued mapping
 $F: E \longrightarrow \mathcal{P}_{cp}(E)$  satisfies the following inequality
 \begin{eqnarray*}\label{disttheorem}
 H_d(Fx,Fy)\leq \varphi(d(x,y)), \ \  \text{     for all    }    x, y     \text {  in }  E,
 \end{eqnarray*}
 where $ \varphi : {\mathbb{R}_+^n} \longrightarrow {\mathbb{R}_+^n} $ is continuous and  nondecreasing function, then
 \begin{eqnarray*}\label{muftheorem}
 \mu(F(V)) \leq  \varphi(\mu(V)), \ \  \text{  for all bounded subset  } V \hbox{  of  }  E.
\end{eqnarray*}
\end{Lemma}
\begin{pf}
 Let $D$ be a bounded subset of $E$ and let  $r :=\left(%
\begin{array}{ccc}
r_1  \\
  \vdots \\
 r_n  \\
\end{array}%
\right),$ with $r_i\geq 0, i = 1,\ldots,n,$ such that $$ \mu(D)=r.$$
By the definition of $\mu,$ for all $i =1,\ldots,n,$ we have
 $$\mu_i(D)= r_i.$$
  Let $ \varepsilon:=\left(%
\begin{array}{ccc}
\varepsilon_1  \\
  \vdots \\
 \varepsilon_n  \\
\end{array}%
\right),$ with $\varepsilon_i> 0, i = 1,\ldots ,n.$
 By the definition of the Hausdorff measure of noncompactness $\mu_i,$ it follows that
 $$
D \subseteq \displaystyle\cup_{k=1}^{m}B_i(x_k,r_i+ \varepsilon_i), \ \  x_k \in E, k=1,\ldots, m,$$
where $B_i(x_k,r_i+ \varepsilon_i)=\left\{ x \in E, \ \ \|x-x_k\|_i \leq r_i+\varepsilon_i  \right\}.$\\
Then, we get
\begin{eqnarray}\label{inclu1}
F(D) \subseteq F\left(\displaystyle\cup_{k=1}^{m}B_i(x_k,r_i+ \varepsilon_i)\right)\subseteq \displaystyle\cup_{k=1}^{m}F \left(B_i(x_k,r_i+ \varepsilon_i)\right).
\end{eqnarray}
Let us consider
$\varphi :=\left(%
\begin{array}{ccc}
\varphi_1  \\
  \vdots \\
 \varphi_n  \\
\end{array}%
\right),$ where
$\varphi_i:\mathbb{R}_+^n \longrightarrow \mathbb{R}_+.$
Now, we claim that
\begin{eqnarray}\label{inclu2}
F\left(B_i(x_k,r_i+ \varepsilon_i)\right) \subseteq  B_i \left(y_k,\varphi_i(r+\varepsilon)\right),\ \  \text{for some } y_k \in \overline{F(x_k)}.
\end{eqnarray}
 Let $y \in F \left(B_i(x_k,r_i+ \varepsilon_i)\right),$ then there exists $x \in B_i(x_k,r_i+\varepsilon_i)$ such that, $y \in Fx.$ Hence, we have
$$\|x-x_k\|_i\leq r_i+\varepsilon_i,$$
which implies that
$$\|x-x_k\| \leq r+\varepsilon.$$
Since $\varphi_i$ is nondecreasing, it follows that
 $$\varphi_i(\|x-x_k\|) \leq \varphi_i( r+\varepsilon).$$
From our assumptions,  we deduce that
\begin{eqnarray*}
H_{d_i}(F(x),F(x_k)) &\leq & \varphi_i(\|x-x_k\| ) \\& \leq &\varphi_i(r+\varepsilon).
\end{eqnarray*}
This means that
\begin{eqnarray*}
d_i(y,F(x_k))& =& \displaystyle\inf_{z_k \in F(x_k)}\|y-z_k\|_i \leq  \displaystyle\inf_{z_k \in \overline{F(x_k)}}\|y-z_k\|_i\\
         &\leq & H_{d_i}(F(x),\overline{F(x_k}))=  H_{d_i}(F(x),F(x_k)) \leq \varphi_i(r+\varepsilon)  ,\  i = 1,\ldots,n.
\end{eqnarray*}
On the other hand,  since ${F(x_k)}$ is compact in $E,$ it follows that
$$\mu((F(x_k)))=0_{\mathbb{R}^n}.$$
Then, $$\mu_i(F(x_k))=0,  \ i = 1,\ldots,n.$$  This implies that $\overline{F(x_k)} $ is compact on Banach space $(E,\|\cdot\|_i),$ and consequently
   there exists $y_k \in \overline{F(x_k)}$   such that
$$ \displaystyle\inf_{z_k \in \overline{F(x_k)}}\|y-z_k\|_i= \|y-y_k\|_i\leq  \varphi_i(r+\varepsilon)$$
and the claim is proved. Now, taking into account both inclusions $(\ref{inclu1})$ and $(\ref{inclu2}),$ we deduce that
 $$F(D)\subseteq  \displaystyle\cup_{k=1}^{m}B_i\left(y_k,\varphi_i(r+\varepsilon) \right),\ i = 1,\ldots,n.$$
Then by again using the definition of the measure of noncompactness $\mu_i,$ we infer that
  $$\mu_i(F(D))\leq \varphi_i(r+\varepsilon),\ i = 1,\ldots,n.  $$
    Since $\varepsilon $ is arbitrary and $\varphi_i$ is continuous, we infer that
$$ \mu_i(F(D)) \leq \varphi_i(\mu(D)),\ i = 1,\ldots,n.$$ Consequently, we get
$$\mu(F(D))=\left(%
\begin{array}{ccc}
\mu_1(F(D))  \\
  \vdots \\
 \mu_n(F(D))  \\
\end{array}%
\right)  \leq \left(%
\begin{array}{ccc}
\varphi_1(\mu(D))  \\
  \vdots \\
\varphi_n( \mu(D))  \\
\end{array}%
\right) = \varphi(\mu(D)).$$
\end{pf}

\begin{remark}
 The condition $\varphi(0)= 0$
 can be removed if we suppose that $Q$ has a closed graph.
\end{remark}
\begin{Theorem}\label{nadler}
Let  $S$ be a nonempty, bounded, closed and convex subset of a generalized Banach space $E.$  Assume that   $F : E \longrightarrow \mathcal{P}_{cp,cv}(E)$ and $G : S \longrightarrow \mathcal{P}(E)$ are two multivalued mappings with closed graphs,  satisfying:\\
\noindent{$(i)$} There is a continuous, nondecreasing function $\varphi:\mathbb{R}^n_+ \longrightarrow  \mathbb{R}^n_+$  with $\lim_{m\longrightarrow\infty}\varphi^m(t)= 0_{\mathbb{R}^n},$ for each $t\in {\mathbb{R}^n_+}$ such that
 \begin{eqnarray*}\label{distanceF}
H_d(Fx,Fy)\leq \varphi(d(x,y)), \ \  \text{     for all    }    x, y     \text {  in }  E,
\end{eqnarray*}
\noindent{$(ii)$} $F(E)$ is bounded and $G(S)$ is relatively compact,\\
\noindent{$(iii)$} $(I-F)^{-1}Gx$ is convex, for all $x \in S,$\\
\noindent {$(iv)$} $x \in F(x)+G(y), y \in S \Longrightarrow  x \in S.$\\
Then $G+F$ has, at least, one fixed point in $S.$\hf
\end{Theorem}
\bigskip
If the multivalued $F$ is $M$-contraction in Nadler's sense, the following result is obtained immediately from Theorem \ref{nadler}.
\begin{Corollary}\label{nadlercontractioncor}
Let  $S$ be a nonempty, bounded, closed and convex subset of a generalized Banach space $X.$  Assume that   $F: E \longrightarrow   \mathcal{P}_{cp, cv}(E)$ and $G: S \longrightarrow   \mathcal{P}(E)$ are multivalued mappings satisfying the following conditions:\\
\noindent{$(i)$}  $F$ is $M$-contraction in  Nadler's sense and $F(E)$ is bounded,\\
\noindent{$(ii)$}  $G$ has a closed graph and  $G(S)$ is relatively compact,\\
\noindent{$(iii)$} $(I-F)^{-1}Gx$ is convex, for all $x \in S,$\\
\noindent {$(iv)$} $x \in F(x)+G(y), y \in S \Rightarrow  x \in S.$\\
Then $G+F$ has, at least, one fixed point in $S.$\hf
\end{Corollary}
\begin{remark}
 Corollary \ref{nadlercontractioncor} extends Theorem $2.14$ in \cite{IOAN}.    \hf
\end{remark}
\bigskip
\section{\textbf{Application to a system of functional-integral inclusions}}
\noindent  Let  $(X,\|\cdot\|)$ be a real Banach space.
  The aim of this section is to discuss existence results
for the following system of functional-integral inclusions in $E \times E:$
\begin{equation}\label{aplicationmulti}
\left\{
  \begin{array}{ll}
    x(t) \in  R_1(t,x(t))+ \displaystyle\int_{0}^{t}Q_1(s,y(s))ds; \\
  \end{array}
\right.
\end{equation}
where $ E:=C(J,X)$ is the space  of all continuous functions  endowed with the norm $\|\cdot\|_\infty= \sup_{t\in J}\|x(t)\|,$ $J:=[0, \rho],$ $\rho>0.$ \\
The following result is important
in our proof:

\noindent In our considerations, we list the following hypotheses:\\
\noindent{($\mathcal{H}_1$)}  The function $R_i:J\times X \times X \longrightarrow X,$ for $i=1,2,$ is such that:\\
\hspace*{20pt}{$($a$)$} There exists a continuous function $h_i:J \longrightarrow \mathbb{R}$ such that
$$ \|R_i(t,x,y)-R_i(s,{x}, {y})\| \leq |h_i(t)-h_i(s)|\, \max(\|x\|, \|y\|), \text{ for all } t,s \in J \text{ and all }x, y \in X.$$
\hspace*{20pt}{$($b$)$} The partial function $(x, y) \mapsto R_i(t,x,y)$ is completely continuous.\\
\noindent{($\mathcal{H}_2$)} The multivalued mapping  $Q_i:J\times X \times X\longrightarrow \mathcal{P}_{cl,cv}(X),$ for $i=1,2,$ is such that:\\
\hspace*{20pt} {$($a$)$}  The partial $Q_{i,t}=Q_i(t, \cdot, \cdot)$  has a closed graph, for all $t\in J.$ \\
\hspace*{20pt} {$($b$)$} For each continuous functions $x, y: J \longrightarrow X,$ there exists a Bochner integrable function
\hspace*{20pt} $ w_i:J \longrightarrow X$ with $w_i(t) \in Q_i(t,x,y)$ a.e on $J.$\\

\begin{Theorem}
Suppose that the assumptions $(\mathcal{H}_{1})-(\mathcal{H}_2)$ are satisfied, the inclusions system (\ref{aplicationmulti}) has at least one solution 
\end{Theorem}

\end{document}